\documentclass[a4paper,12pt]{article}
\usepackage[cp1250]{inputenc}
\usepackage[T1]{fontenc}
\usepackage{amsfonts}
\usepackage{enumerate}
\usepackage{color}
\usepackage{amsmath}
\usepackage{amssymb}
\usepackage{amsthm} 
\usepackage{textcomp}
\usepackage{wrapfig}
\usepackage{url}
\usepackage{fancyhdr}
\usepackage{geometry}
\usepackage{stmaryrd}
\usepackage{bbding}
\usepackage{enumitem} 
\usepackage{bbm}
\everymath{\displaystyle}
\theoremstyle{plain} 
\newtheorem{pr}{Problem}

\theoremstyle{definition} 

\newcommand\cA{{\mathcal A}}
\newcommand\cF{{\mathcal F}}
\newcommand\cS{{\mathcal S}}
\newcommand\cM{{\mathcal M}}

\newcommand\bI{\mathbf{I}_S}

\renewcommand\ge{\geqslant}
\renewcommand\le{\leqslant}

\begin{document}
\begin{center} 
{\LARGE The solution of Hutník's  open problem}
\end{center}

\begin{center}
Michał Boczek
\footnote{Corresponding author. E-mail adress: 800401@edu.p.lodz.pl; 
}, 
Marek Kaluszka 
\\
{\emph{\small{Institute of Mathematics, Lodz University of Technology, 90-924 Lodz, Poland}}}
\end{center}


\begin{abstract} 
In this note, 
we give a solution to Problem $9.2$,  which was presented by Mesiar and Stup\v{n}anov\'{a} [{\it Open problems from the 12th International Conference on Fuzzy Set Theory and Its Applications}, Fuzzy Sets and Systems (2014), http://dx.doi.org/10.1016/ j.fss.2014.07.012]. We show that  
the class of semicopulas solving Problem $9.2$ contains only 
the Łukasiewicz t-norm.
\end{abstract}


\section{The main result}
Let $(X,\cA)$ be a~measurable space, where $\cA$ is a~$\sigma$-algebra of subsets of non-empty set $X,$ and let $\cS$ be the family of all measurable spaces. The~class of all $\cA$-measurable functions $f\colon X\to [0,1]$ is denoted by $\cF_{(X,\cA)}$. A {\it capacity} on $\cA$  is a non-decreasing set function 
$\mu\colon \cA\to [0,1]$ with $\mu(\emptyset)=0$ and  $\mu(X)=1.$
We denote by $\cM_{(X,\cA)}$ the class of all 
capacities  on $\cA.$  Let $S\colon [0,1]^2\to [0,1]$ be a~semicopula (also called 
a $t$-{\it seminorm}), i.e., a~non-decreasing function in both coordinates with the neutral element equal to $1,$ and satisfying the inequality $S(x,y)\le x\wedge y$ for all $x,y\in [0,1],$ where $x\wedge y=\min(x,y)$ (see $\cite{bas}$, $\cite{dur}$ and $\cite{klement2}$). By the above assumptions it follows that $S(x,0)=0=S(0,x)$ for all $x.$
There are three important examples of semicopulas: $M,$ $\Pi$ and $S_L,$ where $M(a,b)=a\wedge b,$ $\Pi(a,b)=ab$ and $S_L(a,b)=(a+b-1)\vee 0;$
$S_L$ is called the {\it Łukasiewicz t-norm}  $\cite{lukasiewicz}.$ Hereafter, 
$a\vee b=\max(a,b)$.

A~class of the smallest semicopula-based integrals is given by
\begin{align*}
 \bI(\mu,f)=\sup_{t\in [0,1]} S\Big(t,\mu\big(\lbrace f\ge t\rbrace \big)\Big),
 \end{align*} 
where $\lbrace f\ge t\rbrace=\lbrace x\in X\colon f(x)\ge t\rbrace,$  $(X,\cA)\in \cS$ and $(\mu, f)\in \cM_{(X,\cA)}\times \cF_{(X,\cA)}.$
Replacing semicopula $S$ with $M$, we get the Sugeno integral $\cite{sugeno1}$. Moreover, if  $S=\Pi,$ then $S$ is called the Shilkret integral $\cite{shilkret}.$

Below we present Problem $9.2$ from $\cite{mesiar11},$ which was posed by Hutník  during the conference FSTA 2014
\textit{The Twelfth International Conference on Fuzzy Set Theory and Applications}
 held from January 26 to January 31, 2014 in Liptovsk\'y J\'an, Slovakia.

 \begin{pr}
$($O. Hutník$)$ Characterize a~class of semicopulas $S,$ for which the equality 
\begin{align}\label{e1}
\bI(\mu,f+a)=\bI(\mu,f)+a
\end{align} 
holds for each $(X,\cA)\in \cS,$ each $(\mu,f)\in \cM_{(X,\cA)}\times \cF_{(X,\cA)}$ and $a\in [0,1]$ such that $f+a\in [0,1].$

\end{pr}

Hal\v{c}inová, Hutník and Molnárová $\cite{hutnik}$ showed that 
the integral equality (1) holds if $S$ is the Łukasiewicz t-norm. 
We show that the class of semicopulas solving Problem 1 contains only 
the Łukasiewicz t-norm.

Firstly, we can observe that 
$\mu\big(\lbrace f+a\ge t\rbrace \big)=1$ for $t\in [0,a],$ so  the left-hand side of $\eqref{e1}$ has the form
\begin{align}
\bI(\mu,f+a)&=\sup_{t\in [0,1]} S\Big(t,\mu\big(\lbrace f+a\ge t\rbrace \big)\Big)\nonumber\\
&=\max\Big[ \sup_{t\in [0,a]} S(t,1), \sup_{t\in (a,1]} S\Big(t,\mu\big(\lbrace f+a\ge t\rbrace \big)\Big)\Big]\nonumber\\&=a\vee \sup_{t\in (0,1-a]} S\Big(t+a,\mu\big(\lbrace f\ge t\rbrace \big)\Big),\label{e2}
\end{align}
since $S(a,1)=a$ for all $a$. Let 
\begin{align}\label{e3}
\mu\big(\lbrace f\ge t\rbrace \big) =
  \begin{cases}
   1 & \text{if } t=0,  \\
   b & \text{if } t\in (0,1-a],  \\
   0      & \text{if } t\in (1-a,1], 
  \end{cases}
\end{align}
where $b\in [0,1].$ By $\eqref{e2}$ and $\eqref{e3}$, 
equality  $\eqref{e1}$ can be rewritten as follows
\begin{align*}
a\vee \sup_{t\in (0,1-a]} S(t+a,b)=\sup_{t\in (0,1-a]}S(t,b) +a.
\end{align*}
As $S$ is a non-decreasing function, we have
\begin{align*}
a\vee S(1,b)= S(1-a,b)+a
\end{align*}
for all $a,b\in [0,1].$ Now, let $c=1-a\in [0,1].$ Thus
 $(1-c)\vee b= S(c,b)+(1-c).$ If $1-c\ge b,$ then
$S(c,b)=0.$ If $1-c< b,$  then $S(c,b)=c+b-1.$  Summing up, $S(c,b)=S_L(c,b)$ for all $c,b\in [0,1].$
\medskip

{\bf Acknowledgments} 
\medskip

This paper has been partially supported
by the grant for young researchers from Lodz University of Technology, grant number $501\backslash 17$-2-2-$7154\backslash 2014.$

\end{document}